\documentclass[11pt,reqno]{amsart}
 
\usepackage{amssymb,latexsym}
\usepackage{graphicx}
\usepackage{url}		 

\calclayout

\makeatletter
\newcommand{\monthyear}[1]{%
  \def\@monthyear{\uppercase{#1}}}
\newcommand{\volnumber}[1]{%
  \def\@volnumber{\uppercase{#1}}}
\AtBeginDocument{%
\def\ps@plain{\ps@empty
  \def\@oddfoot{\@monthyear \hfil \thepage}%
  \def\@evenfoot{\thepage \hfil \@volnumber}}
\def\ps@firstpage{\ps@plain}
\def\ps@headings{\ps@empty
  \def\@evenhead{%
    \setTrue{runhead}%
    \def\thanks{\protect\thanks@warning}%
    \uppercase{Three Fibonacci-chain Aperiodic Algebras}\hfil}%
  \def\@oddhead{%
    \setTrue{runhead}%
    \def\thanks{\protect\thanks@warning}%
    \hfill\uppercase{Three Fibonacci-chain Aperiodic Algebras }}%
  \let\@mkboth\markboth
  \def\@evenfoot{%
    \thepage \hfil \@volnumber}%
  \def\@oddfoot{%
    \@monthyear \hfil \thepage}%
  }%
\footskip=25pt
\pagestyle{headings}%
}
\makeatother

\theoremstyle{plain}
\numberwithin{equation}{section}

\begin{document}
\monthyear{Month Year} \volnumber{Volume, Number} \setcounter{page}{1}
\title{three fibonacci-chain aperiodic algebras }
\author{Daniele Corradetti}
\address{Departimento de Matematica\\
 Universidade do Algarve\\
 Campus de Gambelas,\\
 Faro, PT}
\email{d.corradetti@gmail.com}
\thanks{Research supported by Quantum Gravity Research fundings}
\author{David Chester}
\address{Quantum Gravity Research\\
 Los Angeles, California \\
 CA 90290, USA\\
}
\email{DavidC@QuantumGravityResearch.org}
\author{Raymond Aschheim}
\address{Quantum Gravity Research\\
 Los Angeles, California \\
 CA 90290, USA\\
}
\email{Raymond@QuantumGravityResearch.org}
\author{Klee Irwin}
\address{Quantum Gravity Research\\
 Los Angeles, California \\
 CA 90290, USA\\
}
\email{Klee@QuantumGravityResearch.org}
\begin{abstract}
Aperiodic algebras are infinite dimensional algebras with generators
corresponding to an element of the aperiodic set. These algebras proved
to be an useful tool in studying elementary excitations that can propagate
in multilayered structures and in the construction of some integrable
models in quantum mechanics. Starting from the works of Patera and
Twarock we present three aperiodic algebras based on Fibonacci-chain
quasicrystals: a quasicrystal Lie algebra, an aperiodic Witt algebra
and, finally, an aperiodic Jordan algebra. While a quasicrystal Lie
algebra was already constructed from a modification of the Fibonacci chain,
we here present an aperiodic algebra that matches exactly the original
quasicrystal. Moreover, this is the first time to our knowledge, that
an aperiodic Jordan algebra is presented leaving room for both theoretical
and applicative developments. 
\end{abstract}

\maketitle

\section{Introduction}

Crystallographic Coxeter groups are an essential tool in Lie theory
being in one-to-one correspondence with semisimple finite-dimensional
Lie algebras over the complex number field and, thus, playing a fundamental
role in physics. On the other hand, non-crystallographic Coxeter groups
are deeply connected with icosahedral quasicrystals and numerous aperiodic
structures \cite{MP93,KF15}. It is then natural, in such context,
to look at algebras that are invariant by non-crystallographic symmetries
and, more specifically, at aperiodic Lie algebras. 

Aperiodic algebras are a class of infinite dimensional algebras with
each generator corresponding to an element of the aperiodic set. A
family of aperiodic Lie algebras was firstly introduced by Patera,
Pelantova and Twarock \cite{PPT98}, then generalized and extended
in \cite{PT99,TW00a} and, later on, studied by Mazorchuk \cite{Ma02,MT03}.
These algebras turned out to be suitable for physical applications
and theoretical models such as the breaking of Virasoro symmetry in
\cite{Tw99a} and the construction of exactly solvable models in \cite{TW00b}.
In fact, this is not surprising since important results were obtained
studying elementary excitations that can propagate in multilayered
structures with constituents arranged in a quasiperiodic fashion.
These excitations include plasmon--polaritons, spin waves, light
waves and electrons, among others \cite{Albuque2003}. In this context,
relevant physical properties are analysed in terms of Hamiltonians\cite{Ja21,SCP}
which in the case of nearest-neighbour, tight-binding models are of
the form
\begin{equation}
\left(H\psi\right)_{n}=t\psi_{n+1}+t\psi_{n-1}+\lambda V_{n}\psi_{n},
\end{equation}
where $\lambda$ measures the strength potential, and the potential
sequence $V_{n}$ is generated according to some aperiodic substitution
rule \cite{Macia2006} such as Fibonacci, Thue-Morse, Period-doubling,
Triadic Cantor, etc. It is therefore with renewed interest that we
look to aperiodic algebras that match exactly those chains and especially
the Fibonacci chains.

In this work, we introduce three aperiodic algebras for a class of
Fibonacci-chain quasicrystals\cite{LS86}: the first one is a quasicrystal
Lie algebra, the second is an aperiodic Witt algebra with a Virasoro
extension and, finally, we present an aperiodic Jordan algebra. 

The present work is structured as follows. In sec. 2, we briefly review
the 1-dimensional \emph{Fibonacci-chain} quasicrystals in a general
setting following \cite{LS86}. We then define the binary operation
of \emph{quasiaddition} which encode a geometrical invariance of this
class of quasicrystals. In sec. 3, 4 and 5, we construct three different
Fibonacci-chain aperiodic algebras. The first algebra, treated in
sec. 3, is a quasicrystal Lie algebra and is a review of the original
works of Twarock and Patera\cite{PPT98,TW00a}, which constitutes
our starting point. As for the second aperiodic algebra, i.e. the
aperiodic Witt algebra and its Virasoro extension treated in sec.
4, we used a slightly different approach than \cite{TW00a}. Indeed,
by focusing on the analytic definition of the Fibonacci-chain quasicrystal, we
defined the aperiodic Witt algebra following an index scheme based
on integers rather than aperiodic quasicrystal coordinates. In the
definition of this algebra, we avoided the point defect
introduced by \cite{PPT98,TW00a}, so that the quasicrystal, on which
our aperiodic algebra is based, contains only tiles that have length
$\tau$ and $\tau^{2}$, and matches exactly a Fibonacci-chain quasicrystal.
Finally, in sec. 5 we present the third aperiodic algebra which is
a Jordan algebra realized exploiting the property of invariance of
quasiaddition over Fibonacci-chain quasicrystals. To our knowledge,
this is the first time that an aperiodic Jordan algebra is presented
in literature. 
\begin{figure}
\centering{}\includegraphics[scale=0.4]{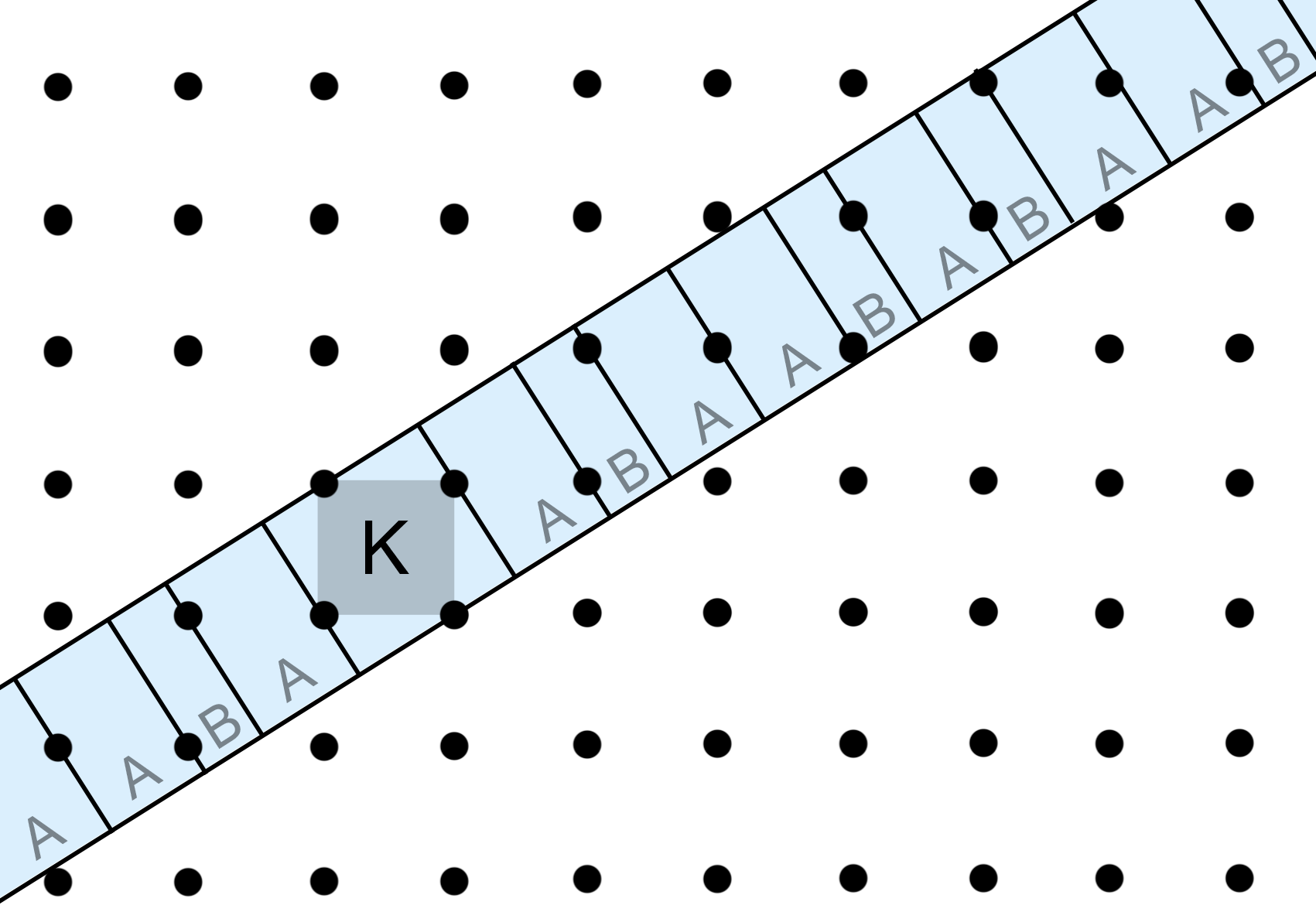}\caption{\label{fig:The-Fibonacci-chain-Quasicrystal}A one-dimensional  quasicrystal
$\mathcal{F}$ obtained through the \emph{cut-and-project scheme}.
The one-dimensional quasicrystal is obtained intersecting an integral
lattice $\mathbb{Z}^{2}$ with the acceptance window $\Omega$ here
represented by a region bounded by two lines of irrational slope $1/\tau$.
Points that fall into the acceptance window are projected onto
the lower line on which lies quasicrystal.}
\end{figure}

\section{The Fibonacci-chain Quasicrystal }

\subsection{Definition}

A \emph{Fibonacci chain} $\mathcal{F}$ is a one-dimensional aperiodic
sequence that can be obtained by the substitution rules $A\longrightarrow AB,\,\,B\longrightarrow A,$
with starting points $S_{0}=B$, which then yields to the following
sequence 
\begin{equation}
B,\,A,\,AB,\,ABA,\,ABAAB,\,...\label{eq: sequenza Fib AB}
\end{equation}
It is easy to note that the number of blocks in (\ref{eq: sequenza Fib AB})
increases following the Fibonacci sequence $1,1,2,3,5...$ and that
the ratio between $A$ and $B$ tends to the \emph{golden mean} $\tau=\left(1+\sqrt{5}\right)/2$.

\subsection{Realisations of the Fibonacci-chain Quasicrystal}

We will now give different realisations of the Fibonacci chain $\mathcal{F}$
in one-dimensional quasicrystals starting from a cut-and-project scheme.
As illustrated in Fig. (\ref{fig:The-Fibonacci-chain-Quasicrystal})
we will consider an integral lattice $\mathbb{Z}^{2}\subset\mathbb{R}^{2}$
and a region, called \emph{acceptance window}, bounded by two lines
with same irrational slope $1/\tau\approx0.618...$ and separated
by an interval that we will suppose here to be of length $1$. Points of
the lattice that fall into the acceptance window will be then projected
to the lower line thus realising a one-dimensional Fibonacci-chain
quasicrystal. Obviously translations of the acceptance window will
correspond to different realisation of the same Fibonacci-chain quasicrystal.

In order to make this construction theoretically precise we here define
a cut-and-project scheme, but we also give an explicit formula for
the coordinates of the quasicrystal. Let $\mathbb{Z}^{2}\subset\mathbb{R}^{2}$
be the integral lattice and $\mathbb{Z}\left[\tau\right]=\mathbb{Z}+\tau\mathbb{Z}$
the Dirichlet ring. Notice that the two roots of the equation $x^{2}=x+1$
are $\tau\approx1.618...$ and $1-\tau\approx-0.618...$, thus
naturally defines a Galois automorphism $*$ that for every number
$n\in\mathbb{Z}\left[\tau\right]$ corresponds
\begin{equation}
n=n_{1}+n_{2}\tau\longrightarrow n^{*}=n_{1}+(1-\tau)n_{2}.\label{eq:starmap}
\end{equation}
The Galois automorphism in (\ref{eq:starmap}) is suitable to be the
\emph{star map} of our cut-and-project quasicrystal, i.e. the involution
that sends points from the parallel space to the perpendicular space
and vice versa\cite{BG17}. We then have that a Fibonacci-chain quasicrystal
is given by

\begin{equation}
\mathcal{F}\left(\Omega\right)=\left\{ n\in\mathbb{Z}\left[\tau\right];n^{*}\in\Omega\right\} ,\label{eq:Model set definition}
\end{equation}
where the acceptance window $\Omega$ is the segment $\left(0,1\right]$.
\begin{table}
\centering{}%
\begin{tabular}{|c|c|c|c|c|c|c|c|c|c|c|c|c|}
\hline 
$n$ &  & $\cdots$ & -4 & -3 & -2 & -1 & 0 & 1 & 2 & 3 & 4 & $\cdots$\tabularnewline
\hline 
$\mathcal{F}_{1,0}\left(n\right)$ &  & $\cdots$ & {\small{}$-2-4\tau$} & {\small{}$-1-3\tau$} & {\small{}$-1-2\tau$} & {\small{}$-\tau$} & {\small{}1} & {\small{}$1+\tau$} & {\small{}$2+2\tau$} & {\small{}$2+3\tau$} & {\small{}$3+4\tau$} & $\cdots$\tabularnewline
\hline 
$\mathcal{F}_{1/2,0}\left(n\right)$ &  & $\cdots$ & {\small{}$-2-4\tau$} & {\small{}$-2-3\tau$} & {\small{}$-1-2\tau$} & {\small{}$-1-\tau$} & {\small{}0} & {\small{}$1+\tau$} & {\small{}$1+2\tau$} & {\small{}$2+3\tau$} & {\small{}$2+4\tau$} & $\cdots$\tabularnewline
\hline 
$\mathcal{F}_{0,0}\left(n\right)$ &  & $\cdots$ & {\small{}$-3-4\tau$} & {\small{}$-2-3\tau$} & {\small{}$-2-2\tau$} & {\small{}$-1-\tau$} & {\small{}0} & {\small{}$\tau$} & {\small{}$1+2\tau$} & {\small{}$1+3\tau$} & {\small{}$2+4\tau$} & $\cdots$\tabularnewline
\hline 
\end{tabular}\caption{\label{tab:Elements-of-the}Elements of the Fibonacci-chain Quasicrystal
$\mathcal{F}_{\alpha,\beta}\left(n\right)$ for $\alpha\in\left\{ 1,\frac{1}{2},0\right\} $
and $\beta=0$.}
\end{table}

It is useful to remark that, by virtue of (\ref{eq:Model set definition}),
a Dirichlet integer $n\in\mathbb{Z}\left[\tau\right]$ belongs to
the Fibonacci-chain quasicrystal if and only if the characteristic
function $\chi_{\Omega}\left(n^{*}\right)=1$. This notation will
come in handy in order to define the Witt aperiodic algebra in the
next section. 

A more explicit way to define the previous Fibonacci-chain 
quasicrystal $\mathcal{F}$ is to consider a direct formulation of the coordinates
as in \cite{LS86}, where a general class of Fibonacci-chain quasicrystals
is defined as the set $\mathcal{F}_{\alpha,\beta}$ containing the
Dirichlet numbers that are image of
\begin{equation}
\mathcal{F}_{\alpha,\beta}\left(m\right)=\left\lfloor \frac{m}{\tau}+\alpha\right\rfloor +m\tau+\beta,\label{eq:Explicit Fab}
\end{equation}
where $m, \beta\in\mathbb{Z}$, $\alpha\in\mathbb{R}$, and $\left\lfloor x\right\rfloor $
is the floor function of $x$. Since the term $\beta$ acts as a translation,
we focus mostly on $\beta=0$ and our primary interest will be in
$\alpha=0,1/2$ or $1$. Indeed, it is straightforward to see that
for $\alpha=1$ and $\beta=0$, then the set is equal to $\mathcal{F}\left(\Omega\right)$
with $\Omega=\left(0,1\right]$ for which a list of explicit values
are given in the first row of Tab. \ref{tab:Elements-of-the}. Another
notable Fibonacci-chain quasicrystal is given by the set $\left\{ \mathcal{F}_{1/2,0}\left(n\right)\right\} _{n\in\mathbb{Z}}$
which is also called the \emph{palindrome Fibonacci-chain} quasicrystal
and whose elements can be also be found in Tab. \ref{tab:Elements-of-the}.
Such quasicrystal has 180-degree rotational symmetry about the origin
and is sometimes convenient for generalizations to higher dimensions.
More generally, a straightforward calculation shows that $\mathcal{F}_{\alpha,\beta}$
identifies the set $\mathcal{F}\left(\Omega\right)$ with $\Omega=\left(-1+\alpha+\beta,\alpha+\beta\right]$.
We thus have that the characteristic function $\chi_{\alpha,\beta}$
for the quasicrystal $\mathcal{F}_{\alpha,\beta}$ is given by 
\begin{equation}
\chi_{\alpha,\beta}\left(n\right)=\left\{ \begin{array}{cc}
1\,\,\,\,\,\, & \mbox{if }-1+\alpha+\beta<n\leq\alpha+\beta\\
0\,\,\,\,\,\, & \mbox{else}
\end{array}\right.,
\end{equation}
where $n\in\mathbb{Z}\left[\tau\right]$.

\subsection{Quasiaddition}

An important symmetry of the Fibonacci-chain quasicrystals is expressed
through a binary operation called \emph{quasiaddition} \cite{MP93} and defined
as

\begin{align}
n & \vdash m=\tau^{2}n-\tau m,
\end{align}
for every $n,m\in\mathbb{Z}\left[\sqrt{5}\right]$. Obviously quasiaddition
is not commutative nor associative but is flexible, i.e. $n\vdash\left(m\vdash n\right)=\left(n\vdash m\right)\vdash n,$and
enjoys the following properties 

\begin{align}
n\vdash n & =n,\\
n\vdash\left(n\vdash m\right) & =m\vdash n,\\
\left(n+p\right)\vdash\left(m+p\right) & =\left(n\vdash m\right)+p,\\
\left(n\vdash m\right)+\left(m\vdash n\right) & =n+m,\\
\left(n\vdash m\right)-\left(m\vdash n\right) & =\left(n-m\right)\vdash\left(m-n\right),
\end{align}
for every $n,m,p\in\mathbb{Z}\left[\sqrt{5}\right]$. An important
feature of this operation is that given any two points $n,m$ of a
Fibonacci-chain quasicrystal $\mathcal{F}\left(\Omega\right)$ then
$n\vdash m$ still belongs to $\mathcal{F}\left(\Omega\right)$. Indeed,
if the image of the star map $n^{*},m^{*}$ in (\eqref{eq:starmap})
belongs to $\Omega$, then also 
\begin{equation}
\left(n\vdash m\right)^{*}=\left(1-\tau\right)^{2}n-\left(1-\tau\right)m,\label{eq:convexF}
\end{equation}
belongs to $\Omega$ since the set is convex. Thus a Fibonacci-chain
quasicrystal is closed under quasiaddition. A sample of the multiplication
table of the quasiaddition is given in Tab. \ref{tab:Quasiaddition-of-}.
\begin{table}
\centering{}{\footnotesize{}}%
\begin{tabular}{|c|c|c|c|c|c|c|c|c|c|}
\hline 
{\footnotesize{}$x\vdash y$} & {\footnotesize{}$y$} & {\footnotesize{}$-1-3\tau$} & {\footnotesize{}$-1-2\tau$} & {\footnotesize{}$-\tau$} & {\footnotesize{}$1$} & {\footnotesize{}$1+\tau$} & {\footnotesize{}$2+2\tau$} & {\footnotesize{}$2+3\tau$} & {\footnotesize{}$\cdots$}\tabularnewline
\hline 
\hline 
{\footnotesize{}$x$} &  & {\footnotesize{}$\vdots$} & {\footnotesize{}$\vdots$} & {\footnotesize{}$\vdots$} & {\footnotesize{}$\vdots$} & {\footnotesize{}$\vdots$} & {\footnotesize{}$\vdots$} & {\footnotesize{}$\vdots$} & \tabularnewline
\hline 
{\footnotesize{}$-1-3\tau$} & {\footnotesize{}$\cdots$} & {\footnotesize{}$-1-3\tau$} & {\footnotesize{}$-2-4\tau$} & {\footnotesize{}$-3-6\tau$} & {\footnotesize{}$-4-8\tau$} & {\footnotesize{}$-5-9\tau$} & {\footnotesize{}$-6-11\tau$} & {\footnotesize{}$-7-12\tau$} & {\footnotesize{}$\cdots$}\tabularnewline
\hline 
{\footnotesize{}$-1-2\tau$} & {\footnotesize{}$\cdots$} & {\footnotesize{}$-\tau$} & {\footnotesize{}$-1-2\tau$} & {\footnotesize{}$-2-4\tau$} & {\footnotesize{}$-3-6\tau$} & {\footnotesize{}$-4-7\tau$} & {\footnotesize{}$-5-9\tau$} & {\footnotesize{}$-6-11\tau$} & {\footnotesize{}$\cdots$}\tabularnewline
\hline 
{\footnotesize{}$-\tau$} & {\footnotesize{}$\cdots$} & {\footnotesize{}$2+2\tau$} & {\footnotesize{}$1+\tau$} & {\footnotesize{}$-\tau$} & {\footnotesize{}$-1-3\tau$} & {\footnotesize{}$-2-4\tau$} & {\footnotesize{}$-3-6\tau$} & {\footnotesize{}$-4-7\tau$} & {\footnotesize{}$\cdots$}\tabularnewline
\hline 
{\footnotesize{}$1$} & {\footnotesize{}$\cdots$} & {\footnotesize{}$4+5\tau$} & {\footnotesize{}$3+4\tau$} & {\footnotesize{}$2+2\tau$} & {\footnotesize{}$1$} & {\footnotesize{}$-\tau$} & {\footnotesize{}$-1-3\tau$} & {\footnotesize{}$-2-4\tau$} & {\footnotesize{}$\cdots$}\tabularnewline
\hline 
{\footnotesize{}$1+\tau$} & {\footnotesize{}$\cdots$} & {\footnotesize{}$5+7\tau$} & {\footnotesize{}$4+6\tau$} & {\footnotesize{}$3+4\tau$} & {\footnotesize{}$2+2\tau$} & {\footnotesize{}$1+\tau$} & {\footnotesize{}$-\tau$} & {\footnotesize{}$-1-2\tau$} & {\footnotesize{}$\cdots$}\tabularnewline
\hline 
{\footnotesize{}$2+2\tau$} & {\footnotesize{}$\cdots$} & {\footnotesize{}$7+11\tau$} & {\footnotesize{}$6+9\tau$} & {\footnotesize{}$5+7\tau$} & {\footnotesize{}$4+5\tau$} & {\footnotesize{}$3+4\tau$} & {\footnotesize{}$2+2\tau$} & {\footnotesize{}$1+\tau$} & {\footnotesize{}$\cdots$}\tabularnewline
\hline 
{\footnotesize{}$2+3\tau$} & {\footnotesize{}$\cdots$} & {\footnotesize{}$8+12\tau$} & {\footnotesize{}$7+11\tau$} & {\footnotesize{}$6+9\tau$} & {\footnotesize{}$5+7\tau$} & {\footnotesize{}$4+6\tau$} & {\footnotesize{}$3+4\tau$} & {\footnotesize{}$2+3\tau$} & {\footnotesize{}$\cdots$}\tabularnewline
\hline 
{\footnotesize{}$\vdots$} &  & {\footnotesize{}$\vdots$} & {\footnotesize{}$\vdots$} & {\footnotesize{}$\vdots$} & {\footnotesize{}$\vdots$} & {\footnotesize{}$\vdots$} & {\footnotesize{}$\vdots$} & {\footnotesize{}$\vdots$} & \tabularnewline
\hline 
\end{tabular}\caption{\label{tab:Quasiaddition-of-}Quasiaddition of $x\vdash y$ where
$x$ and $y$ belong to the Fibonacci-chain Quasicrystal $\mathcal{F}\left(n\right)$.}
\end{table}

\section{A Fibonacci-Chain Quasicrystal Lie Algebras }

In the following sections we present three aperiodic algebras with
generators in one to one correspondence with the Fibonacci-chain quasicrystal
previously defined. The first quasicrystal Lie algebra is just a review
of that presented in \cite{PPT98} for a one-dimensional quasicrystal
originated by a modification of the Fibonacci-chain and that enjoys
a reflection symmetry at $x=1/2$. The tile at the origin is a defect
with length 1, while all other tiles have length $\tau$ and $\tau^{2}$.
The coordinates of this chain are those of $\mathcal{F}_{1,0}\left(n\right)$
(see Tab. \eqref{tab:Elements-of-the}) with the addition of the $0$
element, i.e. 
\begin{table}
\begin{centering}
{\small{}}%
\begin{tabular}{|c||c|c|c|c|c|}
\hline 
{\small{}$\left[L_{x},L_{y}\right]$} & {\small{}$L_{0}$} & {\small{}$L_{1}$} & {\small{}$L_{1+\tau}$} & {\small{}$L_{2+2\tau}$} & {\small{}$L_{2+3\tau}$}\tabularnewline
\hline 
\hline 
{\small{}$L_{-2-4\tau}$} & {\small{}$(2+4\tau)L_{-2-4\tau}$} & {\small{}$0$} & {\small{}$(3+5\tau)L_{-1-3\tau}$} & {\small{}$0$} & {\small{}$(4+7\tau)L_{-\tau}$}\tabularnewline
\hline 
{\small{}$L_{-1-3\tau}$} & {\small{}$(1+3\tau)L_{-1-3\tau}$} & {\small{}$0$} & {\small{}$0$} & {\small{}$0$} & {\small{}$(3+6\tau)L_{1}$}\tabularnewline
\hline 
{\small{}$L_{-1-2\tau}$} & {\small{}$(1+2\tau)L_{-1-2\tau}$} & {\small{}$0$} & {\small{}$(2+3\tau)L_{-\tau}$} & {\small{}$(3+4\tau)L_{1}$} & {\small{}$(3+5\tau)L_{1+\tau}$}\tabularnewline
\hline 
{\small{}$L_{-\tau}$} & {\small{}$\tau L_{-\tau}$} & {\small{}$0$} & {\small{}$(1+2\tau)L_{1}$} & {\small{}$0$} & {\small{}$(2+4\tau)L_{2+2\tau}$}\tabularnewline
\hline 
{\small{}$L_{0}$} & {\small{}$0$} & {\small{}$L_{1}$} & {\small{}$(1+\tau)L_{1+\tau}$} & {\small{}$(2+2\tau)L_{2+2\tau}$} & {\small{}$(2+3\tau)L_{2+3\tau}$}\tabularnewline
\hline 
{\small{}$L_{1}$} & {\small{}$-L_{1}$} & {\small{}$0$} & {\small{}$0$} & {\small{}$0$} & {\small{}$0$}\tabularnewline
\hline 
{\small{}$L_{1+\tau}$} & {\small{}$(-1-\tau)L_{1+\tau}$} & {\small{}$0$} & {\small{}$0$} & {\small{}$0$} & {\small{}$(1+2\tau)L_{3+4\tau}$}\tabularnewline
\hline 
{\small{}$L_{2+2\tau}$} & {\small{}$(-2-2\tau)L_{2+2\tau}$} & {\small{}$0$} & {\small{}$0$} & {\small{}$0$} & {\small{}$\tau L_{4+5\tau}$}\tabularnewline
\hline 
{\small{}$L_{2+3\tau}$} & {\small{}$(-2-3\tau)L_{2+3\tau}$} & {\small{}$0$} & {\small{}$(-1-2\tau)L_{3+4\tau}$} & {\small{}$-\tau L_{4+5\tau}$} & {\small{}$0$}\tabularnewline
\hline 
\end{tabular}\caption{\label{tab:CommL(F(O))TwPat}Commutation relations for a sample of
generators of $\mathfrak{L}\left(\mathcal{F}\left(\Omega\right)\right)$
that are given for various $x,y\in\mathcal{F}\left(\Omega\right)$.}
\par\end{centering}
\label{table:FibChainDefectComm}
\end{table}
 
\begin{equation}
...,-1-3\tau,-1-2\tau,-\tau,0,1,1+\tau,2+2\tau,...
\end{equation}
Following \cite{PPT98}, we now define the quasicrystal Lie algebra
for the Fibonacci-chain quasicrystal $\mathfrak{L}\left(\mathcal{F}\left(\Omega\right)\right)$
as the infinite dimensional vector space spanned by $\left\{ L_{x}\right\} _{x\in\mathcal{F}\left(\Omega\right)}$
with the bilinear product defined by
\begin{equation}
\left[L_{x},L_{y}\right]=\left(y-x\right)\chi_{\Omega}\left(x^{*}+y^{*}\right)L_{x+y},\label{eq:Twarock QCLieAlgebra}
\end{equation}
where $x,y\in\mathcal{F}\left(\Omega\right)$ and $\Omega=\left[0,1\right]$
is the acceptance window of the quasicrystal. Since (\ref{eq:Twarock QCLieAlgebra})
it is obviously antisymmetric, we only have to show that it satisfies
the Jacobi identity. In fact, this is straightforward since 
\begin{align}
\left[L_{x},\left[L_{y},L_{z}\right]\right] & =\left(z-y\right)\left(y+z-x\right)\chi_{\Omega}\left(x^{*}+y^{*}+z^{*}\right)\chi_{\Omega}\left(y^{*}+z^{*}\right)L_{x+y+z},\\
\left[L_{y},\left[L_{z},L_{x}\right]\right] & =\left(x-z\right)\left(z+x-y\right)\chi_{\Omega}\left(x^{*}+y^{*}+z^{*}\right)\chi_{\Omega}\left(z^{*}+x^{*}\right)L_{x+y+z},\\
\left[L_{z},\left[L_{x},L_{y}\right]\right] & =\left(y-x\right)\left(x+y+z\right)\chi_{\Omega}\left(x^{*}+y^{*}+z^{*}\right)\chi_{\Omega}\left(x^{*}+y^{*}\right)L_{x+y+z},
\end{align}
and since $\chi_{\Omega}\left(x^{*}+y^{*}+z^{*}\right)=1$ implies
$\chi_{\Omega}\left(x^{*}+y^{*}\right)=1$ for every $x^{*},y^{*},z^{*}\in\Omega=\left[0,1\right]$
and, more generally for $\Omega=\left[a,b\right]$ with $ab\geq0$.
From previous equations, we have that 
\begin{equation}
\left[L_{x},\left[L_{y},L_{z}\right]\right]+\left[L_{y},\left[L_{z},L_{x}\right]\right]+\left[L_{z},\left[L_{x},L_{y}\right]\right]=0,
\end{equation}
thus fulfilling the Jacobi identity. Therefore, definition (\ref{eq:Twarock QCLieAlgebra})
is a Lie algebra for every $\Omega=\left[a,b\right]$ with $ab\geq0$.
A sample of commutation relations $\left[L_{x},L_{y}\right]$ with
$x,y\in\mathcal{F}\left(\Omega\right)$ are explicitly given in Tab.
\ref{tab:CommL(F(O))TwPat}. 

Observing (\ref{eq:Twarock QCLieAlgebra}), few remarks are easily
spotted. First of all, we notice that if we consider an arbitrary
interval for the acceptance window $\Xi=\left[a,1\right]$, then if
$a\geq1/2$ the Lie algebra is abelian since $\chi_{\Xi}\left(x^{*}+y^{*}\right)$
is tautologically zero for every $x,y\in\mathcal{F}\left(\Xi\right)$.
Moreover, the subalgebra $\mathfrak{L}\left(\mathcal{F}\left(\Xi\right)\right)$
is an ideal of $\mathfrak{L}\left(\mathcal{F}\left(\Omega\right)\right)$
if $\Xi=\left[c,1\right]$ where $0<c$ . This implies that $\mathfrak{L}\left(\mathcal{F}\left(\Omega\right)\right)$
are never semisimple Lie algebras. 
\begin{table}
\centering{}{\small{}}%
\begin{tabular}{|c||c|c|c|c|c|c|c|c|}
\hline 
{\small{}$\left[L_{m},L_{n}\right]$} & {\small{}$L_{0}$} & {\small{}$L_{1}$} & {\small{}$L_{2}$} & {\small{}$L_{3}$} & {\small{}$L_{4}$} & {\small{}$L_{5}$} & {\small{}$L_{6}$} & {\small{}$L_{7}$}\tabularnewline
\hline 
\hline 
{\small{}$L_{-4}$} & {\small{}$-4L_{-4}$} & {\small{}$0$} & {\small{}$-6L_{-2}$} & {\small{}$0$} & {\small{}$0$} & {\small{}$-9L_{1}$} & {\small{}$0$} & {\small{}$-11L_{3}$}\tabularnewline
\hline 
{\small{}$L_{-3}$} & {\small{}$-3L_{-3}$} & {\small{}$-4L_{-2}$} & {\small{}$-5L_{-1}$} & {\small{}$0$} & {\small{}$-7L_{1}$} & {\small{}$-8L_{2}$} & {\small{}$-9L_{3}$} & {\small{}$-10L_{4}$}\tabularnewline
\hline 
{\small{}$L_{-2}$} & {\small{}$-2L_{-2}$} & {\small{}$0$} & {\small{}$0$} & {\small{}$0$} & {\small{}$0$} & {\small{}$-7L_{3}$} & {\small{}$0$} & {\small{}$0$}\tabularnewline
\hline 
{\small{}$L_{-1}$} & {\small{}$-L_{-1}$} & {\small{}$0$} & {\small{}$-3L_{1}$} & {\small{}$0$} & {\small{}$-5L_{3}$} & {\small{}$-6L_{4}$} & {\small{}$0$} & {\small{}$-8L_{6}$}\tabularnewline
\hline 
{\small{}$L_{0}$} & {\small{}$0$} & {\small{}$-L_{1}$} & {\small{}$-2L_{2}$} & {\small{}$-3L_{3}$} & {\small{}$-4L_{4}$} & {\small{}$-5L_{5}$} & {\small{}$-6L_{6}$} & {\small{}$-7L_{7}$}\tabularnewline
\hline 
{\small{}$L_{1}$} & {\small{}$L_{1}$} & {\small{}$0$} & {\small{}$-L_{3}$} & {\small{}$0$} & {\small{}$0$} & {\small{}$-4L_{6}$} & {\small{}$0$} & {\small{}$-6L_{8}$}\tabularnewline
\hline 
{\small{}$L_{2}$} & {\small{}$2L_{2}$} & {\small{}$L_{3}$} & {\small{}$0$} & {\small{}$0$} & {\small{}$-2L_{6}$} & {\small{}$-3L_{7}$} & {\small{}$-4L_{8}$} & {\small{}$-5L_{9}$}\tabularnewline
\hline 
{\small{}$L_{3}$} & {\small{}$3L_{3}$} & {\small{}$0$} & {\small{}$0$} & {\small{}$0$} & {\small{}$0$} & {\small{}$-2L_{8}$} & {\small{}$0$} & {\small{}$0$}\tabularnewline
\hline 
{\small{}$L_{4}$} & {\small{}$4L_{4}$} & {\small{}$0$} & {\small{}$2L_{6}$} & {\small{}$0$} & {\small{}$0$} & {\small{}$-L_{9}$} & {\small{}$0$} & {\small{}$-3L_{11}$}\tabularnewline
\hline 
\end{tabular}\caption{Commutation relations for a sample of generators of the aperiodic
Witt algebras $\mathfrak{W}\left(\mathcal{F}_{0,0}\right)$.}
\label{table:FibChainQCWittComm0}
\end{table}
 
\begin{table}
\centering{}{\small{}}%
\begin{tabular}{|c||c|c|c|c|c|c|c|c|}
\hline 
{\small{}$\left[L_{m},L_{n}\right]$} & {\small{}$L_{0}$} & {\small{}$L_{1}$} & {\small{}$L_{2}$} & {\small{}$L_{3}$} & {\small{}$L_{4}$} & {\small{}$L_{5}$} & {\small{}$L_{6}$} & {\small{}$L_{7}$}\tabularnewline
\hline 
\hline 
{\small{}$L_{-4}$} & {\small{}$0$} & {\small{}$-5L_{-3}$} & {\small{}$0$} & {\small{}$-7L_{-1}$} & {\small{}$-8L_{0}$} & {\small{}$0$} & {\small{}$-10L_{2}$} & {\small{}$0$}\tabularnewline
\hline 
{\small{}$L_{-3}$} & {\small{}$0$} & {\small{}$0$} & {\small{}$0$} & {\small{}$-6L_{0}$} & {\small{}$0$} & {\small{}$0$} & {\small{}$0$} & {\small{}$0$}\tabularnewline
\hline 
{\small{}$L_{-2}$} & {\small{}$0$} & {\small{}$-3L_{-1}$} & {\small{}$-4L_{0}$} & {\small{}$-5L_{1}$} & {\small{}$-6L_{2}$} & {\small{}$0$} & {\small{}$-8L_{4}$} & {\small{}$-9L_{5}$}\tabularnewline
\hline 
{\small{}$L_{-1}$} & {\small{}$0$} & {\small{}$-2L_{0}$} & {\small{}$0$} & {\small{}$-4L_{2}$} & {\small{}$0$} & {\small{}$0$} & {\small{}$-7L_{5}$} & {\small{}$0$}\tabularnewline
\hline 
{\small{}$L_{0}$} & {\small{}$0$} & {\small{}$0$} & {\small{}$0$} & {\small{}$0$} & {\small{}$0$} & {\small{}$0$} & {\small{}$0$} & {\small{}$0$}\tabularnewline
\hline 
{\small{}$L_{1}$} & {\small{}$0$} & {\small{}$0$} & {\small{}$0$} & {\small{}$-2L_{4}$} & {\small{}$-3L_{5}$} & {\small{}$0$} & {\small{}$-5L_{7}$} & {\small{}$0$}\tabularnewline
\hline 
{\small{}$L_{2}$} & {\small{}$0$} & {\small{}$0$} & {\small{}$0$} & {\small{}$-L_{5}$} & {\small{}$0$} & {\small{}$0$} & {\small{}$0$} & {\small{}$0$}\tabularnewline
\hline 
{\small{}$L_{3}$} & {\small{}$0$} & {\small{}$2L_{4}$} & {\small{}$L_{5}$} & {\small{}$0$} & {\small{}$-L_{7}$} & {\small{}$0$} & {\small{}$-3L_{9}$} & {\small{}$-4L_{10}$}\tabularnewline
\hline 
{\small{}$L_{4}$} & {\small{}$0$} & {\small{}$3L_{5}$} & {\small{}$0$} & {\small{}$L_{7}$} & {\small{}$0$} & {\small{}$0$} & {\small{}$-2L_{10}$} & {\small{}$0$}\tabularnewline
\hline 
\end{tabular}\caption{Commutation relations for a sample of generators of the aperiodic
Witt algebras $\mathfrak{W}\left(\mathcal{F}_{1,0}\right)$.}
\label{table:FibChainQCWittComm1}
\end{table}

\section{An aperiodic Witt algebra and its Virasoro extensions }

We now define an aperiodic Witt algebra and a Virasoro extension of
it. Such algebras were first introduced in \cite{TW00a}, but are
here presented in an equivalent way over the base $\left\{ L_{n}\right\} _{n\in\mathbb{Z}}$
with integer index, instead of over the base $\left\{ L_{x}\right\} _{x\in\mathcal{F}_{\alpha,\beta}}$
which was indexed by points of the Fibonacci-chain quasicrystal. First
of all we have to notice from (\ref{eq:Explicit Fab}) that ${F}_{\alpha,\beta}$ is of the form
\begin{equation}
\mathcal{F}_{\alpha,\beta}\left(n\right)=n'+\tau n,
\end{equation}
for some $n'\in\mathbb{Z}$. Thus, the correspondence between $\left\{ L_{n}\right\} _{n\in\mathbb{Z}}$
and $\left\{ L_{x}\right\} _{x\in\mathcal{F}_{\alpha,\beta}}$ is
easily obtained by setting $\mathcal{F}_{\alpha,\beta}\left(n\right)\longrightarrow n\in\mathbb{Z}$.
Therefore we define the aperiodic Witt algebra $\mathfrak{W}\left(\mathcal{F}_{\alpha,\beta}\right)$
as the vector space spanned by $\left\{ L_{n}\right\} _{n\in\mathbb{Z}}$,
equipped with the following bilinear product
\begin{equation}
\left[L_{n},L_{m}\right]=\left(n-m\right)\chi_{\Omega}\left(\mathcal{F}_{\alpha,\beta}\left(n\right)^{*}+\mathcal{F}_{\alpha,\beta}\left(m\right)^{*}\right)L_{n+m}.\label{eq:aperWitt}
\end{equation}
A straightforward calculation, similar to the one presented in the
previous section, shows that this is indeed a Lie algebra if $\left(-1+\alpha+\beta\right)\left(\alpha+\beta\right)\geq0$.
This means that in the special case of the palindrome Fibonacci-chain,
where $\alpha=1/2$ and $\beta=0$, we do not have an aperiodic Witt
algebra. On the other hand, leaving $\beta=0$, we have that $\alpha=0$
and $\alpha=1$ give rise to two different aperiodic Witt algebras
that are non-abelian and, therefore, interesting cases. Starting with
$\alpha=0$ and $\beta=0$, explicit commutation relations for $\mathfrak{W}\left(\mathcal{F}_{0,0}\right)$
are shown in Table \eqref{table:FibChainQCWittComm0} while those
for $\alpha=1$ and $\beta=0$, i.e. $\mathfrak{W}\left(\mathcal{F}_{1,0}\right)$,
are shown in Table \eqref{table:FibChainQCWittComm1}. The structure
constants here are found to be in terms of Dirichlet integers, which
is in general irrational. 

\subsection{An aperiodic Virasoro extension}

We will now focus on the Fibonacci-chain quasicrystals $\mathcal{F}_{\alpha,0}$
with $\alpha=0,1$. While Fibonacci-chain quasicrystal Lie algebras
such as $\mathfrak{L}\left(\mathcal{F}\left(\Omega\right)\right)$
do not allow for a central extension (cfr. \cite{PT99}), the previously
defined aperiodic Witt algebra $\mathfrak{W}\left(\mathcal{F}\left(\Omega\right)\right)$
does. To achieve that central extension in our case, i.e. an aperiodic
Virasoro algebra for the Fibonacci-chain $\mathfrak{V}\left(\mathcal{F}_{\alpha,0}\left(\Omega\right)\right)$,
we just follow \cite{TW00a} adapting it to our special case in our
notation. We then define $\mathfrak{V}\left(\mathcal{F}_{\alpha,0}\left(\Omega\right)\right)$
as the vector space spanned by $C\cup\left\{ L_{n}\right\} _{n\in\mathbb{Z}}$,
where $C$ is the central generator of the extension, with a bilinear
product given by $\left[L_{n},C\right]=0$ for all $L_{n}\in\left\{ L_{n}\right\} _{n\in\mathbb{Z}}$
and
\begin{equation}
\begin{array}{cc}
\left[L_{n},L_{m}\right]=\left(n-m\right) & \chi_{\Omega}\left(\mathcal{F}_{\alpha,0}\left(n\right)^{*}+\mathcal{F}_{\alpha,0}\left(m\right)^{*}\right)L_{n+m}+\frac{1}{12}n\left(n^{2}-1\right)\delta_{n,-m}C,\end{array}
\end{equation}
for every $n,m\in\mathbb{Z}$. Such extension might not be unique,
but is an interesting algebra for its possible applications in theoretical
physics (see a similar algebra in \cite{Tw99a}) and it was then worth
mentioning it. 
\begin{table}
\centering{}{\small{}}%
\begin{tabular}{|c||c|c|c|c|c|c|c|c|}
\hline 
{\small{}$\left[L_{m},L_{n}\right]$} & {\small{}$L_{0}$} & {\small{}$L_{1}$} & {\small{}$L_{2}$} & {\small{}$L_{3}$} & {\small{}$L_{4}$} & {\small{}$L_{5}$} & {\small{}$L_{6}$} & {\small{}$L_{7}$}\tabularnewline
\hline 
\hline 
{\small{}$L_{-4}$} & {\small{}$-4L_{-4}$} & {\small{}$0$} & {\small{}$-6L_{-2}$} & {\small{}$0$} & {\small{}$5C$} & {\small{}$-9L_{1}$} & {\small{}$0$} & {\small{}$-11L_{3}$}\tabularnewline
\hline 
{\small{}$L_{-3}$} & {\small{}$-3L_{-3}$} & {\small{}$-4L_{-2}$} & {\small{}$-5L_{-1}$} & {\small{}$2C$} & {\small{}$-7L_{1}$} & {\small{}$-8L_{2}$} & {\small{}$-9L_{3}$} & {\small{}$-10L_{4}$}\tabularnewline
\hline 
{\small{}$L_{-2}$} & {\small{}$-2L_{-2}$} & {\small{}$0$} & {\small{}$\frac{1}{2}C$} & {\small{}$0$} & {\small{}$0$} & {\small{}$-7L_{3}$} & {\small{}$0$} & {\small{}$0$}\tabularnewline
\hline 
{\small{}$L_{-1}$} & {\small{}$-L_{-1}$} & {\small{}$0$} & {\small{}$-3L_{1}$} & {\small{}$0$} & {\small{}$-5L_{3}$} & {\small{}$-6L_{4}$} & {\small{}$0$} & {\small{}$-8L_{6}$}\tabularnewline
\hline 
{\small{}$L_{0}$} & {\small{}$0$} & {\small{}$-L_{1}$} & {\small{}$-2L_{2}$} & {\small{}$-3L_{3}$} & {\small{}$-4L_{4}$} & {\small{}$-5L_{5}$} & {\small{}$-6L_{6}$} & {\small{}$-7L_{7}$}\tabularnewline
\hline 
{\small{}$L_{1}$} & {\small{}$L_{1}$} & {\small{}$0$} & {\small{}$-L_{3}$} & {\small{}$0$} & {\small{}$0$} & {\small{}$-4L_{6}$} & {\small{}$0$} & {\small{}$-6L_{8}$}\tabularnewline
\hline 
{\small{}$L_{2}$} & {\small{}$2L_{2}$} & {\small{}$L_{3}$} & {\small{}$0$} & {\small{}$0$} & {\small{}$-2L_{6}$} & {\small{}$-3L_{7}$} & {\small{}$-4L_{8}$} & {\small{}$-5L_{9}$}\tabularnewline
\hline 
{\small{}$L_{3}$} & {\small{}$3L_{3}$} & {\small{}$0$} & {\small{}$0$} & {\small{}$0$} & {\small{}$0$} & {\small{}$-2L_{8}$} & {\small{}$0$} & {\small{}$0$}\tabularnewline
\hline 
{\small{}$L_{4}$} & {\small{}$4L_{4}$} & {\small{}$0$} & {\small{}$2L_{6}$} & {\small{}$0$} & {\small{}$0$} & {\small{}$-L_{9}$} & {\small{}$0$} & {\small{}$-3L_{11}$}\tabularnewline
\hline 
\end{tabular}\caption{Commutation relations for a sample of generators of the aperiodic
Virasoro algebra $\mathfrak{V}\left(\mathcal{F}_{0,0}\left(\Omega\right)\right)$.}
\label{table:FibChainQCViraComm0}
\end{table}
 
\begin{table}
\centering{}%
\begin{tabular}{|c||c|c|c|c|c|c|c|c|}
\hline 
{\small{}$\left[L_{m},L_{n}\right]$} & $L_{0}$ & $L_{1}$ & $L_{2}$ & $L_{3}$ & $L_{4}$ & $L_{5}$ & $L_{6}$ & $L_{7}$\tabularnewline
\hline 
\hline 
$L_{-4}$ & $0$ & $-5L_{-3}$ & $0$ & $-7L_{-1}$ & $5C-8L_{0}$ & $0$ & $-10L_{2}$ & $0$\tabularnewline
\hline 
$L_{-3}$ & $0$ & $0$ & $0$ & $2C-6L_{0}$ & $0$ & $0$ & $0$ & $0$\tabularnewline
\hline 
$L_{-2}$ & $0$ & $-3L_{-1}$ & $\frac{1}{2}C-4L_{0}$ & $-5L_{1}$ & $-6L_{2}$ & $0$ & $-8L_{4}$ & $-9L_{5}$\tabularnewline
\hline 
$L_{-1}$ & $0$ & $-2L_{0}$ & $0$ & $-4L_{2}$ & $0$ & $0$ & $-7L_{5}$ & $0$\tabularnewline
\hline 
$L_{0}$ & $0$ & $0$ & $0$ & $0$ & $0$ & $0$ & $0$ & $0$\tabularnewline
\hline 
$L_{1}$ & $0$ & $0$ & $0$ & $-2L_{4}$ & $-3L_{5}$ & $0$ & $-5L_{7}$ & $0$\tabularnewline
\hline 
$L_{2}$ & $0$ & $0$ & $0$ & $-L_{5}$ & $0$ & $0$ & $0$ & $0$\tabularnewline
\hline 
$L_{3}$ & $0$ & $2L_{4}$ & $L_{5}$ & $0$ & $-L_{7}$ & $0$ & $-3L_{9}$ & $-4L_{10}$\tabularnewline
\hline 
$L_{4}$ & $0$ & $3L_{5}$ & $0$ & $L_{7}$ & $0$ & $0$ & $-2L_{10}$ & $0$\tabularnewline
\hline 
\end{tabular}\caption{Commutation relations for a sample of generators of the aperiodic
Virasoro algebra $\mathfrak{V}\left(\mathcal{F}_{1,0}\left(\Omega\right)\right)$.}
\label{table:FibChainQCViraComm1}
\end{table}

\section{A Fibonacci-chain Jordan Algebra }

In this section we define for the first time an aperiodic Jordan algebra,
which is an infinite dimensional Jordan algebra $\mathfrak{J}\left(\mathcal{F}\right)$,
whose generators $\left\{ L_{x}\right\} _{x\in\mathcal{F}}$ are in
one to one correspondence with the Fibonacci-chain quasicrystal $\mathcal{F}$.
Since our construction is valid for all $\alpha$ and $\beta$
we will drop such notation and write just $\mathcal{F}$ unless is
needed. We then define $\mathfrak{J}\left(\mathcal{F}\right)$ as
the real vector space spanned by $\left\{ L_{x}\right\} _{x\in\mathcal{F}\left(\Omega\right)}$
with bilinear product given by

\begin{equation}
L_{x}\circ L_{y}=\frac{1}{2}\left(L_{x\vdash y}+L_{y\vdash x}\right).\label{eq:definizione Jordan}
\end{equation}
The product in (\ref{eq:definizione Jordan}) is obviously well-defined
since for every $x,y\in\mathcal{F}$ then $x\vdash y,y\vdash x\in\mathcal{F}$
as shown in (\ref{eq:convexF}). Moreover, it is clearly commutative,
i.e.
\begin{equation}
L_{x}\circ L_{y}=L_{y}\circ L_{x},
\end{equation}
so, in order to prove that is a Jordan algebra it is sufficient to show
that it satisfies the Jordan identity. Indeed, since $x\vdash x=x$,
we have that 
\begin{align}
\left(L_{x}\circ L_{y}\right)\circ\left(L_{x}\circ L_{x}\right) & =\left(L_{x}\circ L_{y}\right)\circ L_{x}\\
 & =\frac{1}{2}\left(L_{x\vdash y}\circ L_{x}+L_{y\vdash x}\circ L_{x}\right)\\
 & =\frac{1}{2}\left(L_{x}\circ L_{x\vdash y}+L_{x}\circ L_{y\vdash x}\right)\\
 & =L_{x}\circ\left(L_{y}\circ L_{x}\right)\\
 & =L_{x}\circ\left(L_{y}\circ\left(L_{x}\circ L_{x}\right)\right),
\end{align}
so that the Jordan identity 
\begin{equation}
\left(L_{x}\circ L_{y}\right)\circ\left(L_{x}\circ L_{x}\right)=L_{x}\circ\left(L_{y}\circ\left(L_{x}\circ L_{x}\right)\right),
\end{equation}
is fulfilled. Jordan algebras are notorious for an abundance of idempotent
elements and this makes no exception. Indeed, we note that since $x\vdash x=x$
all elements of the basis $\left\{ L_{x}\right\} _{x\in\mathcal{F}_{1/2}}$
are idempotent elements. 

An alternative definition of the algebra can be given on a basis whose
elements are indexed by integer. Indeed, let $\mathcal{F}\left(n\right)=n'+\tau n$
and $\mathcal{F}\left(m\right)=m'+\tau m$, and then consider the
quasiaddition of two elements, i.e.
\begin{equation}
\mathcal{F}\left(n\right)\vdash\mathcal{F}\left(m\right)=\tau^{2}\left(n'+\tau n\right)-\tau\left(m'+\tau m\right).
\end{equation}
A straightforward calculation shows that
\begin{equation}
\mathcal{F}\left(n\right)\vdash\mathcal{F}\left(m\right)=\mathcal{F}\left(n'-m'+2n-m\right),
\end{equation}
while, on the other hand, switching the two addends we obtain
\begin{equation}
\mathcal{F}\left(m\right)\vdash\mathcal{F}\left(n\right)=\mathcal{F}\left(m'-n'+2m-n\right).
\end{equation}
We now have an equivalent definition of $\mathfrak{J}\left(\mathcal{F}\right)$
on the vector space spanned by $\left\{ L_{n}\right\} _{n\in\mathbb{Z}}$
with bilinear product defined as
\begin{equation}
L_{n}\circ L_{m}=\frac{1}{2}\left(L_{n'-m'+2n-m}+L_{m'-n'+2m-n}\right),\label{eq:Jordan product nm}
\end{equation}
where $n',m'\in\mathbb{Z}$ and $n'=\tau n-\mathcal{F}\left(n\right)$
and $m'=\tau m-\mathcal{F}\left(m\right)$. 

A simple analysis of (\ref{eq:Jordan product nm}) shows an useful
property of the above Jordanian product. Let $n,m\in\mathbb{Z}$ and
call $p,q\in\mathbb{Z}$ the indices of the generators such that
\begin{equation}
L_{n}\circ L_{m}=\frac{1}{2}\left(L_{p}+L_{q}\right).
\end{equation}
We thus have from (\ref{eq:Jordan product nm}) and a straighforward
calculation that
\begin{equation}
p+q=n+m.\label{eq:x+y=00003Dn+m}
\end{equation}

In order to give concreteness to our construction we will now focus
on a specific aperiodic Jordan algebra, i.e. $\mathfrak{J}\left(\mathcal{F}_{1,0}\right)$.
For such algebra a sample of its multiplication table is in Tab. \ref{tab:Mult-table-F10-1}.
\begin{table}
\begin{centering}
{\small{}}%
\begin{tabular}{|c|c|c|c|c|c|c|}
\hline 
{\small{}$L_{a}\circ L_{b}$} &  & {\small{}$L_{-2}$} & {\small{}$L_{-1}$} & {\small{}$L_{0}$} & {\small{}$L_{1}$} & {\small{}$L_{2}$}\tabularnewline
\hline 
\hline 
{\small{}$L_{-4}$} &  & {\small{}$\frac{1}{2}\left(L_{1}+L_{-7}\right)$} & {\small{}$\frac{1}{2}\left(L_{4}+L_{-9}\right)$} & {\small{}$\frac{1}{2}\left(L_{7}+L_{-11}\right)$} & {\small{}$\frac{1}{2}\left(L_{9}+L_{-12}\right)$} & {\small{}$\frac{1}{2}\left(L_{12}+L_{-14}\right)$}\tabularnewline
\hline 
{\small{}$L_{-3}$} &  & {\small{}$\frac{1}{2}\left(L_{-1}+L_{-4}\right)$} & {\small{}$\frac{1}{2}\left(L_{2}+L_{-6}\right)$} & {\small{}$\frac{1}{2}\left(L_{5}+L_{-8}\right)$} & {\small{}$\frac{1}{2}\left(L_{7}+L_{-9}\right)$} & {\small{}$\frac{1}{2}\left(L_{10}+L_{-11}\right)$}\tabularnewline
\hline 
{\small{}$L_{-2}$} &  & {\small{}$L_{-2}$} & {\small{}$\frac{1}{2}\left(L_{1}+L_{-4}\right)$} & {\small{}$\frac{1}{2}\left(L_{4}+L_{-6}\right)$} & {\small{}$\frac{1}{2}\left(L_{6}+L_{-7}\right)$} & {\small{}$\frac{1}{2}\left(L_{9}+L_{-9}\right)$}\tabularnewline
\hline 
{\small{}$L_{-1}$} &  & {\small{}$\frac{1}{2}\left(L_{1}+L_{-4}\right)$} & {\small{}$L_{-1}$} & {\small{}$\frac{1}{2}\left(L_{-3}+L_{2}\right)$} & {\small{}$\frac{1}{2}\left(L_{4}+L_{-4}\right)$} & {\small{}$\frac{1}{2}\left(L_{7}+L_{-6}\right)$}\tabularnewline
\hline 
{\small{}$L_{0}$} &  & {\small{}$\frac{1}{2}\left(L_{4}+L_{-6}\right)$} & {\small{}$\frac{1}{2}\left(L_{-3}+L_{2}\right)$} & {\small{}$L_{0}$} & {\small{}$\frac{1}{2}\left(L_{2}+L_{-1}\right)$} & {\small{}$\frac{1}{2}\left(L_{5}+L_{-3}\right)$}\tabularnewline
\hline 
{\small{}$L_{1}$} &  & {\small{}$\frac{1}{2}\left(L_{6}+L_{-7}\right)$} & {\small{}$\frac{1}{2}\left(L_{4}+L_{-4}\right)$} & {\small{}$\frac{1}{2}\left(L_{2}+L_{-1}\right)$} & {\small{}$L_{1}$} & {\small{}$\frac{1}{2}\left(L_{4}+L_{-1}\right)$}\tabularnewline
\hline 
{\small{}$L_{2}$} &  & {\small{}$\frac{1}{2}\left(L_{9}+L_{-9}\right)$} & {\small{}$\frac{1}{2}\left(L_{7}+L_{-6}\right)$} & {\small{}$\frac{1}{2}\left(L_{5}+L_{-3}\right)$} & {\small{}$\frac{1}{2}\left(L_{4}+L_{-1}\right)$} & {\small{}$L_{2}$}\tabularnewline
\hline 
{\small{}$L_{3}$} &  & {\small{}$\frac{1}{2}\left(L_{11}+L_{-10}\right)$} & {\small{}$\frac{1}{2}\left(L_{9}+L_{-7}\right)$} & {\small{}$\frac{1}{2}\left(L_{7}+L_{-4}\right)$} & {\small{}$\frac{1}{2}\left(L_{6}+L_{-2}\right)$} & {\small{}$\frac{1}{2}\left(L_{4}+L_{1}\right)$}\tabularnewline
\hline 
{\small{}$L_{4}$} &  & {\small{}$\frac{1}{2}\left(L_{14}+L_{-12}\right)$} & {\small{}$\frac{1}{2}\left(L_{12}+L_{-9}\right)$} & {\small{}$\frac{1}{2}\left(L_{10}+L_{-6}\right)$} & {\small{}$\frac{1}{2}\left(L_{9}+L_{-4}\right)$} & {\small{}$\frac{1}{2}\left(L_{7}+L_{-1}\right)$}\tabularnewline
\hline 
\end{tabular}{\small\par}
\par\end{centering}
\centering{}\caption{\label{tab:Mult-table-F10-1}Multiplication table for a sample of
generators of the aperiodic Jordan algebra $\mathfrak{J}\left(\mathcal{F}_{1,0}\right)$.}
\end{table}
 An implication of (\ref{eq:x+y=00003Dn+m}), together with the idempotency
of the generators, and the multiplication relations in Tab. \ref{tab:Mult-table-F10-1},
is that the algebra $\mathfrak{J}\left(\mathcal{F}_{1,0}\right)$
is non-unital. Indeed, suppose it exists an identity element $I$.
This would imply that $I\circ L_{0}=L_{0}\circ I=L_{0}$. We would
then have
\begin{align}
L_{0}\circ L_{n_{1}}+L_{0}\circ L_{n_{2}}+...+L_{0}\circ L_{n_{m}} & +...=L_{0}
\end{align}
which means that at least one element $L_{r}$ would give $L_{0}\circ L_{r}=L_{0}$.
But from (\ref{eq:x+y=00003Dn+m}), the only possibility is for $L_{0}=L_{r}$.
Repeating the argument for all $L_{n}$ we have that the identity
element must be of the form $I=...+L_{-1}+L_{0}+L_{1}+L_{2}+...$
. But, if so, then consider 
\begin{equation}
L_{0}\circ\left(...+L_{-1}+L_{0}+L_{1}+L_{2}+...\right).
\end{equation}
Since (\ref{eq:Jordan product nm}) the functions $0\vdash n$ and
$n\vdash0$ are monotone functions in $n$ then $L_{-1}$ resulting
from $L_{0}\circ L_{1}$ (see Tab.\eqref{tab:Mult-table-F10-1}) cannot
be obtained by any other $L_{0}\circ L_{n}$ with $n\in\mathbb{Z}$. Therefore, the $L_{-1}$ component is non zero and, thus, in contraddiction with $I$ being the identity element. 

A similar argument shows that the ideal $\mathfrak{k}=\left\{ L_{0}\circ x:x\in\mathfrak{J}\left(\mathcal{F}_{1,0}\right)\right\} $
is a proper ideal of $\mathfrak{J}\left(\mathcal{F}_{1,0}\right)$
since $L_{1}\notin\mathfrak{k}$. Thus $\mathfrak{J}\left(\mathcal{F}_{1,0}\right)$
is not a simple algebra.

\section{Conclusions and Developments}

In this paper we presented three aperiodic algebras that originate
quite naturally from a special class of Fibonacci-chain quasicrystals.
While the first one is a review of the one introduced in \cite{PPT98}
and is a quasicrystal Lie algebra obtained with a one-point defect
from the Fibonacci-chain quasicrystal $\mathcal{F}_{1,0}$, the second
one is an aperiodic Witt algebra, i.e. $\mathfrak{W}\left(\mathcal{F}_{1,0}\right)$,
that exactly matches such quasicrystal and that can be extended to
a Virasoro algebra $\mathfrak{V}\left(\mathcal{F}_{1,0}\right)$.
Finally, we introduced a completely new class of aperiodic algebras,
i.e. the aperiodic Jordan algebras, and presented a special case for
the same Fibonacci-chain quasicrystal, i.e. $\mathfrak{J}\left(\mathcal{F}_{1,0}\right)$.
Such aperiodic Jordan algebras were made possible by exploiting an
important symmetry of Fibonacci-chain quasicrystals which is encoded
by quasiaddition and holds also for higher dimensional quasicrystals.
The definition of such algebras is already an interesting subject
from a mathematical point of view, but from the physical side we think
the three class of algebras we presented here are an indispensable
tools to everyone who is interested in physical aperiodic structures
that can be modeled by Fibonacci-chain quasicrystals. 

Indeed, aperiodic Virasoro algebras can be used to study deformations of exactly solvable models of Calogero-Sutherland type \cite{TW00b}. More specifically, perturbations of the Hamiltonian describing a many-body quantum mechanical system on a circle with $n$ identical particles of mass $m$ can be expressed in terms of the generators of an aperiodic Virasoro algebra. On the other hand, as for the aperiodic Jordan algebra, it is well known a one-to-one correspondence between finite dimensional Jordan algebras and multifield Korteweg-de Vries equations. Indeed, in the early '90, Svinolupov \cite{Sv91} and Sokolov showed that a generalisation of the KdV equation, i.e.

\begin{equation} u_{t}^{i}=u_{xxx}^{i}-6a_{jk}^{i}u^{j}u_{x}^{k},\label{eq:equazione integrable} \end{equation} 
possesses nondegenerate generalised symmetries or conservation laws if and only if $\left\{ a_{jk}^{i}\right\} $ are constants of structure of a finite dimensional Jordan algebra\cite[Thm 1.1]{Sv93}. Similar results hold for the modified KdV equation\cite{Sv93}, for the Sine-Gordon equation and for generalisations of the non-linear Schroedinger equation\cite{Sv92}. While, the aperiodic Jordan algebra $\mathfrak{J}\left(\mathcal{F}_{1,0}\right)$ is infinite dimensional, nevertheless it can be used to define finite dimensional algebras limiting the aperiodic set of the quasicrystal to a specific set consecutive elements and imposing the product to be $0$ when the quasiaddition of two elements falls outside the selected portion of the quasicrystal. With this modification, commutativity and the Jordan identity still hold, so that the resulting algebra is finite dimensional, falling into Svinolupov theorem's hypothesis. In fact, we expect that systems of the type (\ref{eq:equazione integrable}) where the constants of structure $\left\{ a_{jk}^{i}\right\} $ are those of $\mathfrak{J}\left(\mathcal{F}_{1,0}\right)$ might arise in the study of solitons propagating onto parallel lines departing from points of a $\mathcal{F}_{1,0}$ Fibonacci-chain quasicrystal. 

\section{Acknowledgments}

This work was funded by the Quantum Gravity Research institute. Authors would like to thank Fang Fang, Marcelo Amaral, Dugan Hammock, and Richard Clawson for insightful discussions and suggestions.

\medskip{}

\noindent MSC2020: 11B39, 52C23, 17B65, 17C50 

\begin{thebibliography}{PPT98}
\bibitem[AC03]{Albuque2003}E.L. Albuquerque and M.G. Cottam, \emph{Theory
of elementary excitations in quasiperiodic structures,} Phys. Rep
\textbf{376.4--5} (2003):225-337.

\bibitem[BG17]{BG17}M. Baake and U. Grimm, \emph{Aperiodic Order},
volume 2 of\emph{ Encyclopedia of Mathematics and its Applications},
Cambridge University Press, 2017.

\bibitem[KF15]{KF15}Fang F. and Klee I., \emph{An Icosahedral Quasicrystal
as a Golden Modification of the Icosagrid and its Connection to the
E8 Lattice,} Acta Crystallogr. A: Found. Adv. (2015) 71(a1).

\bibitem[Ja21]{Ja21}A. Jagannathan, \emph{The Fibonacci quasicrystal:
Case study of hidden dimensions and multifractality}, Reviews of Modern
Physics, 93.4 (2021):045001.

\bibitem[Ko99]{Ko99}Koecher M., \emph{The Minnesota Notes on Jordan
Algebras and Their Applications}, Springer Verlag, 1999.

\bibitem[LS86]{LS86}D. Levine and P. J. Steinhardt, \emph{Quasicrystals.
I. definition and structure}, Phys. Rev. B, \textbf{34} (1986):596--616.

\bibitem[Ma]{Macia2006}E. Macia, \emph{The role of aperiodic order
in science and technology}, Rep. Prog. Phys. \textbf{69} (2006):397.

\bibitem[Ma02]{Ma02}V. Mazorchuk, \emph{On quasicrystal lie algebras}.
Journal of Mathematical Physics, 43.5 (2002):2791.

\bibitem[MT03]{MT03}V. Mazorchuk and R. Twarock, \emph{Virasoro-type
algebras associated with a penrose tiling}, Journal of Physics A:
Mathematical and General, \textbf{36.15} (2003):4363-4373.

\bibitem[MP93]{MP93}Moody R. V. and Patera J. \emph{Quasicrystals
and icosians}. Journal of Physics A: Mathematical and General , \textbf{26}
(12):2829, 1993.

\bibitem[PPT98]{PPT98}J. Patera, E. Pelantova and R. Twarock, \emph{Quasicrystal
Lie algebras}, Physics Letters A, \textbf{246.3-4} (1998):209-213. 

\bibitem[PT99]{PT99}J. Patera and R. Twarock, \emph{Quasicrystal
Lie algebras and their generalizations}, Physics Reports, \textbf{315.1-3}
(1999):241-256. 

\bibitem[SCP]{SCP}A. Sen and C. Castro Perelman, \emph{A Hamiltonian
model of the Fibonacci quasicrystal using non-local interactions:
simulations and spectral analysis}, The European Physical Journal
B, \textbf{93.4} (2020). 

\bibitem[Se95]{Sench Quasicryst}M. Senechal, \emph{Quasicrystals
and Geometry}, Cambridge University Press, 1995.

\bibitem[Sv91]{Sv91}Svinolupov S.I., \emph{Jordan algebras and generalized
Korteweg-de Vries equations}. Theor Math Phys \textbf{87} (1991) 611--620.

\bibitem[Sv92]{Sv92}Svinolupov S.I., \emph{Generalized Schrödinger
equations and Jordan pairs}. Commun.Math. Phys. \textbf{143 }(1992)
559--575.

\bibitem[Sv93]{Sv93}Svinolupov S.I., \emph{Jordan algebras and integrable
systems} Funct. Anal. its Appl. \textbf{27} (1993) 257--265

\bibitem[TW00a]{TW00a}R. Twarock, \emph{Aperiodic Virasoro algebra}.
Journal of Mathematical Physics, \textbf{41.7} (2000):5088-5106. 

\bibitem[TW00b]{TW00b}R. Twarock. \emph{An exactly solvable Calogero
model for a non-integrally laced group}. Physics Letters A, \textbf{275.3}
(2000):169172.

\bibitem[Tw99a]{Tw99a}R. Twarock \emph{Breaking Virasoro symmetry
}in \emph{ Quantum Theory and Symmetries} ed H-D Doebner, V K Dobrev,
J-D Hennig and W Lucke, Singapore: World Scientific, 1999, pp 530--4.
\end{thebibliography}
\end{document}